\documentclass[a4paper,DIV=11,abstract=true]{scrartcl}

\usepackage[T1]{fontenc}
\usepackage[english]{babel}
\usepackage{amsmath}
\usepackage{amsfonts}
\usepackage{amssymb}
\usepackage{graphicx}
\usepackage{booktabs}
\usepackage{mathtools}
\usepackage{amsthm}
\usepackage{xcolor}
\usepackage[bf,normal]{caption}
\usepackage{subcaption}
\usepackage[ruled,vlined]{algorithm2e}
\usepackage{enumitem}
\usepackage{authblk}

\definecolor{deeppink}{rgb}{1.0, 0.08, 0.58}
\newcommand{\new}[1]{\textcolor{black}{#1}}

\usepackage[pdffitwindow=false,
			plainpages=false,
			pdfpagelabels=true,
			pdfpagemode=UseOutlines, % how to open PDF externally
			pdfpagelayout=OneColumn, 
			bookmarks=true,
			bookmarksopen=true,
			colorlinks=true,
			hyperfootnotes=false,
			linkcolor=blue,
			urlcolor=blue!30!black,
			citecolor=green!50!black]{hyperref}

\theoremstyle{plain} % italic
\newtheorem{thm}{Theorem}[section]

\theoremstyle{definition} % not italic
\newtheorem{defn}[thm]{Definition}

\newtheorem{rem}[thm]{Remark}
\newtheorem{expl}[thm]{Example}

\newcommand{\ts}{\hspace*{0.1em}} 								            % Thin space
\newcommand{\ddt}[1]{\frac{\mathrm{d}}{\mathrm{d}#1}}                       % Derivatives
\newcommand{\pd}[2]{\frac{\partial#1}{\partial#2}}                          % Partial derivatives
\newcommand{\underbracedmatrix}[2]{%
	\left[\;
	\smash[b]{\underbrace{
			\begin{matrix}#1\end{matrix}
		}_{#2}}
	\;\right]
	\vphantom{\underbrace{\begin{matrix}#1\end{matrix}}_{#2}}
}

\newcommand\xqed[1]{\leavevmode\unskip\penalty9999 \hbox{}\nobreak\hfill \quad\hbox{#1}}
\newcommand{\exampleSymbol}{\xqed{$\triangle$}}

\title{Koopman-Based Surrogate Models \\ for Multi-Objective Optimization \\ of Agent-Based Systems}

\author[1, 2]{Jan-Hendrik Niemann}
\author[3]{Stefan Klus}
\author[1]{\\ Nata\v{s}a Djurdjevac Conrad}
\author[1,2]{Christof Schütte}
\affil[1]{Department of Mathematics and Computer Science, Freie Universit\"at Berlin, Germany}
\affil[2]{Zuse Institute Berlin, Germany}
\affil[3]{School of Mathematical \& Computer Sciences, Heriot--Watt University, UK}

\date{}

\begin{document}

\maketitle

\begin{abstract}
Agent-based models (ABMs) provide an intuitive and powerful framework for studying social dynamics by modeling the interactions of individuals from the perspective of each individual. In addition to simulating and forecasting the dynamics of ABMs, the demand to solve optimization problems to support, for example, decision-making processes naturally arises. Most ABMs, however, are non-deterministic, high-dimensional dynamical systems, so objectives defined in terms of their behavior are computationally expensive. In particular, if the number of agents is large, evaluating the objective functions often becomes prohibitively time-consuming. We consider data-driven reduced models based on the Koopman generator to enable the efficient solution of multi-objective optimization problems involving ABMs. In a first step, we show how to obtain data-driven reduced models of non-deterministic dynamical systems (such as ABMs) that depend potentially \new{nonlinearly on} control inputs. We then use them in the second step as surrogate models to solve multi-objective optimal control problems. We first illustrate our approach using the example of a voter model, where we compute optimal controls to steer the agents to a predetermined majority, and then using the example of an epidemic ABM, where we compute optimal containment strategies in a prototypical situation. We demonstrate that the surrogate models effectively approximate the Pareto-optimal points of the ABM dynamics by comparing the surrogate-based results with test points, where the objectives are evaluated using the ABM. Our results show that when objectives are defined by the dynamic behavior of ABMs, data-driven surrogate models support or even enable the solution of multi-objective optimization problems.
\end{abstract}


\section{Introduction}

Modeling social dynamics and studying the resulting collective phenomena is an important research problem, for instance, in the field of epidemic modeling, opinion dynamics, mobility, or innovation spreading in ancient times \cite{Goldenbogen2022,Banisch2012,Mielke2018,DjurdjevacConrad2018}. Agent-based models (ABMs) provide an intuitive yet powerful framework for modeling the interactions of individuals, small groups, or entire populations. The high level of modeling flexibility allows users to gain important insights into the complex dynamic patterns that emerge from the interactions of discrete entities called \emph{agents}. These agents often follow simple rules that describe their behavior, making ABMs accessible to both experts and non-experts without extensive mathematical knowledge. Agent-based modeling is particularly useful in modeling complex systems where interactions and environmental factors play a significant role, where heterogeneity among agents is desired or required, or where agents can learn and adapt to new situations. In addition to simulating and forecasting the dynamics of complex ABMs, optimization problems such as optimal control need to be solved. One important use case is to support decision-making processes. 

Optimal control is the process of determining the best set of actions or inputs to a dynamical system over a given time horizon in order to achieve a desired objective while satisfying given constraints. There exist well-established numerical methods to solve complex, nonlinear optimization problems such as line search, conjugate gradient, or trust region methods~\cite{NocedalWright2006,Conn2000}. The optimization of ABMs, however, leads to significant challenges since basic concepts such as derivatives are not well-defined due to the fact that most ABMs are inherently stochastic and discontinuous. Additionally, most ABMs are high-dimensional systems comprising thousands of agents, which is often not only desirable but also necessary, e.g., to be able to correctly represent even small fractions of the total population. Furthermore, many independent simulations are often required to compute statistical properties such as means and variances. This renders the simulation of ABMs expensive and time-consuming and a more thorough analysis including optimization and also sensitivity analysis or uncertainty quantification nearly impossible.

In this work, we will focus on so-called \emph{multi-objective optimization problems} involving ABMs. These problems arise naturally whenever multiple objectives are to be optimized simultaneously, without any predetermined prioritization. There exist various methods and techniques to solve multi-objective optimization problems, such as scalarization, $\varepsilon$-constraint methods, evolutionary algorithms, particle swarm, agent-based and multi-agent methods, or set-oriented methods, see~\cite{Miettinen1998, Schuetze2005, CoelloCoello2006, Peitz2018, Afsar2020, Blondin2020, Dellnitz2005, Schuetze2013}.

In the setting of ABMs, one way to solve multi-objective optimization problems is to use heuristic methods such as evolutionary algorithms. These methods do not rely on derivatives and are simple to use, see, e.g.,~\cite{An2017} and references therein, and are commonly used for an automatic or guided calibration of ABMs, for example, to fit the model against real-world financial market data~\cite{Rogers2004} or biological observations~\cite{Calvez2006,Read2016}. Multi-objective calibration can also provide insights into parameters or initial conditions that may not be evident through simulations alone, see~\cite{Liu2017, Moya2019, ParragaAlava2019, KorkmazTan2019, DAuria2020, Bartkowski2020} for recent studies focusing on multi-objective calibration of ABMs using heuristic methods. We refer the reader to~\cite{Platt2020} for a comprehensive review of calibration techniques. However, the lack of (mathematical) convergence properties and the often costly simulation requirements are drawbacks of heuristic methods for solving multi-objective optimization problems involving ABMs.

Surrogates, either for individual objectives or for the complete model, provide a way to reduce the computational effort when optimizing ABMs. In this sense, a \emph{surrogate} refers to a viable substitute \new{for the original model} that has a sufficient level of accuracy and can be evaluated significantly faster, often by several orders of magnitude. Interpolation, regression, or machine learning can be used to quickly obtain a surrogate for expensive objectives~\cite{Peitz2018, Berkemeier2021}. However, when model parameters are changed, these surrogates need to be recalculated. Instead of replacing expensive objective functions, the entire model itself can be replaced by a surrogate model. In~\cite{Lamperti2018}, a surrogate model that approximates the mapping between ABM inputs and the corresponding response in the output is learned using non-parametric machine learning for simulation, calibration, and exploration of the parameter space. Further approaches include, e.g., surrogate models based on difference equations~\cite{Oremland2014, Oremland2015, KoshyChenthittayil2021}, partial differential equations~\cite{Christley2017}, or ordinary differential equations (ODEs)~\cite{Wulkow_H2021} to solve multi-objective optimization problems involving ABMs. See also~\cite{Peitz2018} for a recent review on surrogate modeling.

In this work, we use a data-driven method to find suitable surrogate models. Data-driven methods prove to be particularly advantageous for this task as they can extract valuable insights about the behavior of dynamical systems solely from data, without relying on prior knowledge about the system. This makes these methods particularly well-suited to address problems where a closed description of the system is not available. \new{Various methods have been developed in recent years, e.g., for approximating transfer operators associated with the system~\cite{Klus2016}, for detecting metastable and coherent sets~\cite{schutte_klus_hartmann_2023}, or for performing stability analyses~\cite{Mauroy2016a}, as well as for deriving the governing equations of the underlying dynamics \cite{Mauroy2016, Kutz2016, Klus2020}, model reduction, optimization, and control \cite{Korda2018, Arbabi2018, Korda2020, Peitz2019, Peitz2020}.} Following~\cite{Niemann2021a}, we use the infinitesimal generator of the Koopman operator, which is the adjoint of the Perron--Frobenius operator associated with the dynamical system and describes the evolution of observable functions representing any kind of measurement, to obtain a reduced model from ABM simulation data, which is then used as a surrogate model to solve a multi-objective optimization problem. \new{The reduction of the numerical effort is achieved in two ways: First, a data-driven method allows us to obtain a dynamical system that models only the required quantities. In our particular case these are the aggregated dynamics of an ABM, i.e., the collective behavior of the agents. This corresponds to a dimension reduction. Second, the reduced dynamical model, which serves as a surrogate model, is given in terms of differential equations, which can in general be simulated more efficiently than ABMs.} The main contributions of this work are:
\begin{itemize}[itemsep=0ex, topsep=0.5ex]
    \item We show how \emph{generator extended dynamic mode decomposition} (gEDMD)~\cite{Klus2020}, a data-driven method that approximates the infinitesimal generator of the Koopman operator, can be extended to identify the governing equations of non-deterministic dynamical systems with potentially nonlinear \new{dependence on the} control inputs.
    \item We demonstrate how these reduced models, obtained from ABM simulation data, can be used as surrogate models to enable the efficient solution of multi-objective optimal control problems.
    \item We consider two different use cases, namely linear and nonlinear \new{dependence on the} control inputs for two ABMs, and show for both ABMs that the Pareto sets computed using the data-driven surrogate models indeed approximate the Pareto sets of the ABMs.
\end{itemize}
In this study, we showcase the potential of data-driven surrogate models to efficiently address the challenges posed by multi-objective optimization problems associated with the dynamical behavior of ABMs. In particular when limit models for large numbers of agents are unknown or non-existent, data-driven surrogate models offer a viable solution for tackling optimization problems that would otherwise be computationally infeasible due to the high computational cost of evaluating the objectives. Our method requires a significant amount of data, but in simulation studies where a surrogate model is needed for optimization or control, data availability is typically not an issue.

The remainder of this paper is organized as follows. In Section~\ref{sec:preliminaries} we introduce the stochastic Koopman operator and its generator as well as multi-objective optimization problems. We then show in Section~\ref{sec:Koopman_Reduced_Models_with_Control} how gEDMD can be used to obtain data-driven reduced models with linear and nonlinear \new{dependence on the} control inputs. We demonstrate in Section~\ref{sec:MOOP4ABMS} how both approaches can be used to find surrogate models to solve multi-objective optimal control problems for high-dimensional ABMs. We consider the voter model and an ABM modeling the transmission dynamics of SARS-CoV-2. Open questions and future work will be discussed in Section~\ref{sec:conclusion}.

%%%%%%%%%%%%%%%%%%%%%%%%%%%%%%%%%%%%%%%%%%%%%%%%%%%%%%%%%%%%%%%%%%%%%%%%%%%

\section{Preliminaries}
\label{sec:preliminaries}

We will first introduce the required mathematical tools and provide a brief overview of the stochastic Koopman operator and its generator as well as multi-objective optimization problems.

\subsection{The Koopman Operator and its Generator}
\label{sec:Koopman_operator_generator}

Let $\mathbb{X} \subset \mathbb{R}^d$ be the state space. We consider stochastic differential equations (SDEs) with (time-varying) control input $u \colon \mathbb{R}_{\ge 0} \to \mathbb{R}^{d_u}$ of the form
\begin{equation} \label{eq:SDE+control}
    \mathrm{d}X_t = b(X_t, u) \ts \mathrm{d}t + \sigma(X_t, u) \ts \mathrm{d}W_t,
\end{equation}
where $\{X_t\}_{t \ge 0} \in \mathbb{X}$ is a time-homogeneous stochastic process and $ b \colon \mathbb{R}^d \times \mathbb{R}^{d_u} \to \mathbb{R}^d $ denotes the drift term, $ \sigma \colon \mathbb{R}^d \times \mathbb{R}^{d_u} \to \mathbb{R}^{d \times s} $ the diffusion term, and $ W_t $ an $ s $-dimensional Wiener process. Let $ \Phi^t $ denote the associated flow map \new{and $ f \in L^{\infty}(\mathbb{X}) $ a real-valued observable of the system representing any kind of measurement}. Assuming that $u$ is constant, the semigroup $ \{\ts \mathcal{K}^t_u \ts\}_{t \ge 0} $ of Koopman operators $\mathcal{K}^t_u \colon L^{\infty}(\mathbb{X}) \to L^{\infty}(\mathbb{X})$ is defined by
\begin{equation*}
	(\mathcal{K}^t_u f)(x) = \mathbb{E}[f(\Phi^t(x,u)) \mid X_t = x],
\end{equation*}
\new{i.e., the conditional expectation of $f(\Phi^t(x,u))$ given $X_t = x$}. The Koopman operator is an infinite-dimensional, linear and non-expansive operator \new{forming a contraction semigroup}. Provided that $f$ is twice continuously differentiable, it can be shown using It\^{o}'s lemma that the infinitesimal generator $\mathcal{L}_u$ of the Koopman operator can be characterized by
\begin{equation*}
	(\mathcal{L}_u f)(x)
	= \sum_{i=1}^d b_i(x, u) \ts \pd{f}{x_i} + \frac{1}{2} \sum_{i=1}^d \sum_{j=1}^d a_{ij}(x, u) \ts \pd{^2 f}{x_i \ts \partial x_j},
\end{equation*}
where $a = \sigma \ts \sigma^\top$. We refer to it as the \emph{Koopman generator}\new{, see~\cite{Lasota1994} for details}. The function $ v(t, x) \coloneqq (\mathcal{K}^t_u f)(x) $ satisfies $ \pd{v}{t} = \mathcal{L}_u v $, which is a second-order partial differential equation commonly known as the \emph{Kolmogorov backward equation}, see~\cite{Metzner2007}. For deterministic dynamical systems, i.e., $\sigma \equiv 0 $, we obtain a first-order partial differential equation, the Liouville equation. See~\cite{Klus2020} for details.

\begin{rem}
    We denote the Koopman operator and its generator for an uncontrolled system of the form~\eqref{eq:SDE+control}, i.e., $u \equiv 0$, by $\mathcal{K}_0^t$ and $\mathcal{L}_0$, respectively.
\end{rem}

\subsubsection{Generator Extended Dynamic Mode Decomposition}

Due to the infinite-dimensional nature of the Koopman operator, it is common practice to consider projections onto finite-dimensional subspaces. We briefly introduce \emph{generator extended dynamic mode decomposition} (gEDMD), a method which approximates the Koopman generator. For details, see~\cite{Klus2020}. We will assume for now that $u \equiv 0$ and omit any subscripts to simplify the notation. Given $m$ measurements of the system's state $ \{\ts x_l \ts\}_{l=1}^m $, its drift $ \{\ts b(x_l) \ts\}_{l=1}^m $, and diffusion $ \{\ts \sigma(x_l) \ts\}_{l=1}^m $, as well as a set of basis functions $ \{\ts \psi_k \ts\}_{k=1}^\ell $, which can be written in vector form as $ \psi(x) = [\psi_1(x), \dots, \psi_\ell(x)]^\top $, we define
\begin{equation*}
    \mathrm{d}\psi_k(x) = (\mathcal{L} \psi_k)(x) = \sum_{i=1}^d b_i(x) \ts \pd{\psi_k}{x_i}(x) + \frac{1}{2} \sum_{i=1}^d \sum_{j=1}^d a_{ij}(x) \ts \pd{^2 \psi_k}{x_i \ts \partial x_j}(x)
\end{equation*}
\new{for $k = 1, \dots, \ell $}. We then construct the matrices $ \Psi_X, \mathrm{d}\Psi_X \in \mathbb{R}^{\ell \times m} $ for all measurements and basis functions\new{, i.e.,}
\begin{equation} \label{eq:gEDMDmatrices}
	\Psi_X =
	\begin{bmatrix}
	\psi_1(x_1) & \dots  & \psi_1(x_m) \\
	\vdots      & \ddots & \vdots      \\
	\psi_\ell(x_1) & \dots  & \psi_\ell(x_m)
	\end{bmatrix}
	\quad \text{and} \quad
	\mathrm{d}\Psi_X =
	\begin{bmatrix}
	\mathrm{d}\psi_1(x_1) & \dots  & \mathrm{d}\psi_1(x_m) \\
	\vdots                & \ddots & \vdots                \\
	\mathrm{d}\psi_\ell(x_1) & \dots  & \mathrm{d}\psi_\ell(x_m)
	\end{bmatrix},
\end{equation}
\new{and} assume that $ \mathrm{d}\Psi_X = M \Psi_X $. In general, this problem cannot be solved exactly such that we solve it in the least-square sense by minimizing $ \Vert\mathrm{d}\Psi_X - M \Psi_X\Vert_F $, where $\Vert\cdot\Vert_F$ denotes the Frobenius norm. The solution is given by $M = \mathrm{d}\Psi_X \Psi_X^+$, where $\Psi_X^+$ denotes the Moore--Penrose pseudoinverse of $\Psi_X$. The matrix $ L = M^\top $ is an empirical estimate of the matrix representation of the infinitesimal generator $\mathcal{L}$. The convergence of gEDMD in the infinite data limit to a Galerkin approximation of the generator, i.e., a projection onto the space spanned by the basis functions, was shown in~\cite{Klus2020}.

\subsubsection{System Identification}

Assuming that $\mathbb{X}$ is bounded so that the full-state observable $g(x) = x$ is contained com\-po\-nent-wise in $L^{\infty}(\mathbb{X})$, we can use the observable to reconstruct the governing equations of the underlying dynamical system. We make the assumption that $ g(x) = x $ can be represented by the basis functions $ \psi $, which can easily be accomplished by adding the observables $ \{ \ts x_i \ts \}_{i=1}^d $ to the set of basis functions. The system can directly be represented in terms of the basis functions,
\begin{equation*}
	(\mathcal{L} g)(x) = b(x) \approx (L B)^\top \ts \psi(x),
\end{equation*}	
where \new{we define} $ B \in \mathbb{R}^{\ell \times d} $ \new{such} that $ g(x) = B^\top \ts \psi(x) $. For deterministic dynamical systems, this is equivalent to SINDy~\cite{Brunton2015}. For non-deterministic systems and for $\psi_k(x) = x_ix_j$, we identify the diffusion term using
\begin{equation*}
	a_{ij}(x) \approx (\mathcal{L}\psi_k)(x) - b_i(x) x_j - b_j(x) x_i,
\end{equation*}
where we assume that $b_i$ and $b_j$ as well as $b_i(x)x_j$ and $b_j(x)x_i$ are contained in the space spanned by the basis functions. If it is necessary to obtain the drift term $\sigma$, it can be obtained by computing a Cholesky decomposition of $a$. For details, see~\cite{Klus2020}.

\subsection{Multi-Objective Optimization}
\label{sec:MOOP}

\emph{Multi-objective optimization} concerns the simultaneous optimization of $k$ -- often contradictory -- objective functions $f_1, \dots, f_k \colon \mathbb{R}^n \to \mathbb{R}$. For example, one might want to minimize the spread of an infectious disease with curfews while at the same time also trying to minimize the socioeconomic cost of this intervention. A decision maker must find a compromise, i.e., an agreement reached through mutual concessions that often alter or combine the original goals. An optimal compromise cannot be further improved without worsening at least one of the other goals. This is called \emph{Pareto-optimal}. In contrast to single-objective optimization, where conflicting objectives are prioritized a priori, e.g., by a weighted sum of the objectives, multi-objective optimization prioritizes a posteriori.  This means that if we know Pareto-optimal points, we can select a solution according to the preference of each objective.

We consider multi-objective optimization problems of the form
\begin{equation} \label{eq:MOOP}
	\begin{aligned}
		& \underset{y \in \mathbb{R}^n}{\text{min}} &  F(y) &= \left[f_1(y), \dots, f_k(y)\right]^\top, \\
		& \text{s.t.} &  g_i(y) &\le 0, \quad i=1, \dots, q,\\
		& &  h_j(y) &= 0, \quad j=1, \dots, p,
	\end{aligned}
\end{equation}
where $F \colon \mathbb{R}^n \to \mathbb{R}^k$, $g \colon \mathbb{R}^n \to \mathbb{R}^q$ and $h \colon \mathbb{R}^n \to \mathbb{R}^p$ for $g \coloneqq [g_1, \dots, g_q]^\top$ and $h \coloneqq [g_1, \dots, g_p]^\top$. The space of all feasible decisions is called the \emph{feasible decision space} $\mathcal{R}$ and is given by the constraints in~\eqref{eq:MOOP}, i.e.,
\begin{equation*}
	\mathcal{R} = \{ y \in \mathbb{R}^n \mid g(y) \le 0, h(y) = 0 \},
\end{equation*}
where $\le$ is defined component-wise. A decision $y^\star \in \mathcal{R}$ is called \emph{\new{(globally)} Pareto-optimal} or \emph{\new{(global)} Pareto point} if there is no $y \in \mathcal{R}$ that dominates $y^\star$, i.e., $F(y) \neq F(y^\star)$ and $ F(y) \le F(y^\star)$ hold component-wise. \new{A decision $y^\star$ is called \emph{locally Pareto-optimal} or \emph{local Pareto point} if there exists $\epsilon > 0$ such that $y^\star$ is Pareto-optimal in a neighborhood $y \in U_\epsilon(y^\star) \subset \mathcal{R}$.} The set of all Pareto-optimal points is called \emph{Pareto set} and its image under the map of the objective functions \emph{Pareto front}~\cite{Schuetze2005, Schuetze2013}.

\begin{figure}
	\centering
	\hspace{-0.5cm}
	\begin{subfigure}[t]{0.49\textwidth}
		\centering
		\includegraphics[scale=1]{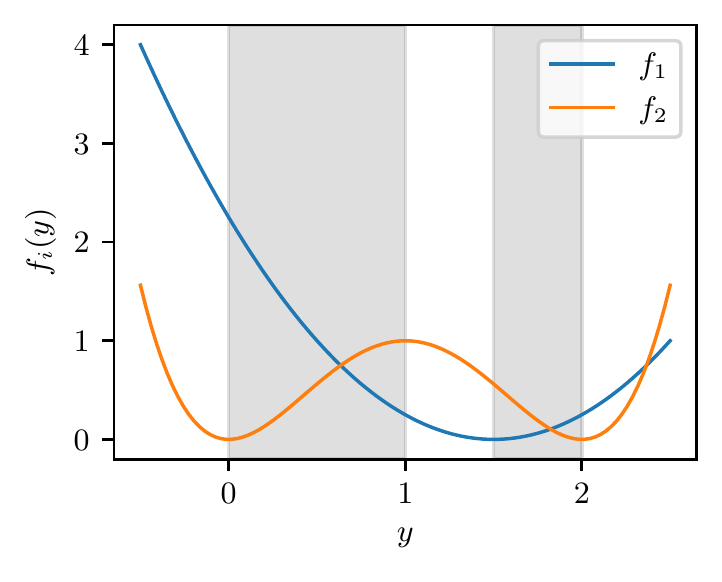}
		\caption{}
	\end{subfigure}
	\hfil
	\begin{subfigure}[t]{0.49\textwidth}
		\centering
		\includegraphics[scale=1]{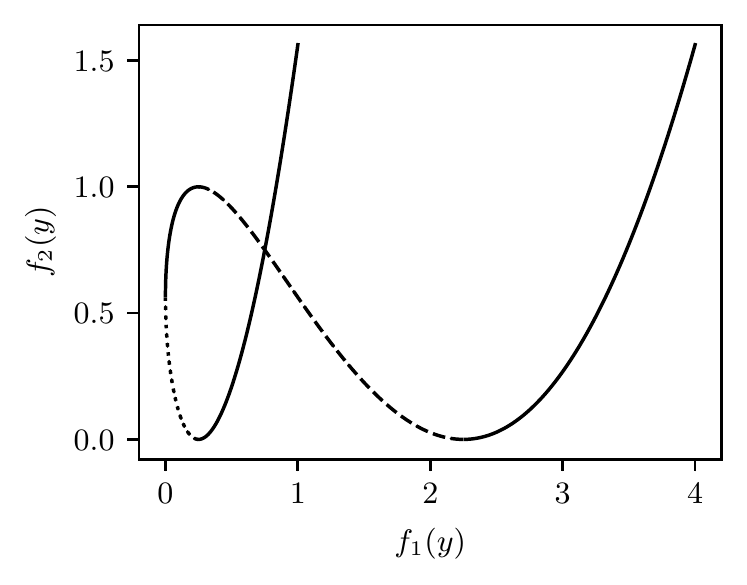}
		\caption{}
	\end{subfigure}
	\caption{(a)~Objective functions $f_1$ and $f_2$ on the interval $\mathcal{R} = [-0.5, 2.5]$ with Pareto set (gray shaded) and (b)~Pareto front (dotted and dashed). The dotted line segment shows that only points in $[1.5, 2]$ are globally Pareto-optimal since they clearly dominate all points in $[0, 1]$.}
	\label{fig:example_sampling_algo}
\end{figure}

\begin{expl} \label{expl:pareto}
	Consider the multi-objective optimization problem~\eqref{eq:MOOP} from~\cite{Dellnitz2005} with objective functions $f_i \colon \mathcal{R} \to \mathbb{R}$ given by
	\begin{align*}
		f_1(y) &= (y - 1.5)^2, \\
		f_2(y) &= y^4 - 4 \ts y^3 + 4 \ts y^2
	\end{align*}
	and feasible decision space $\mathcal{R} = [-0.5, 2.5]$. Figure~\ref{fig:example_sampling_algo}(a) shows both objective functions. The set \new{of locally Pareto-optimal points} is given by $[0, 1] \cup [1.5, 2]$. Figure~\ref{fig:example_sampling_algo}(b) shows not only the Pareto front as dotted and dashed line segments, but also that only the interval $[1.5, 2]$ is globally Pareto-optimal, as every point in $[1.5, 2]$ clearly dominates all points in the other interval. \exampleSymbol
\end{expl}

\begin{figure}[t]
    \centering
    \includegraphics[width=\textwidth]{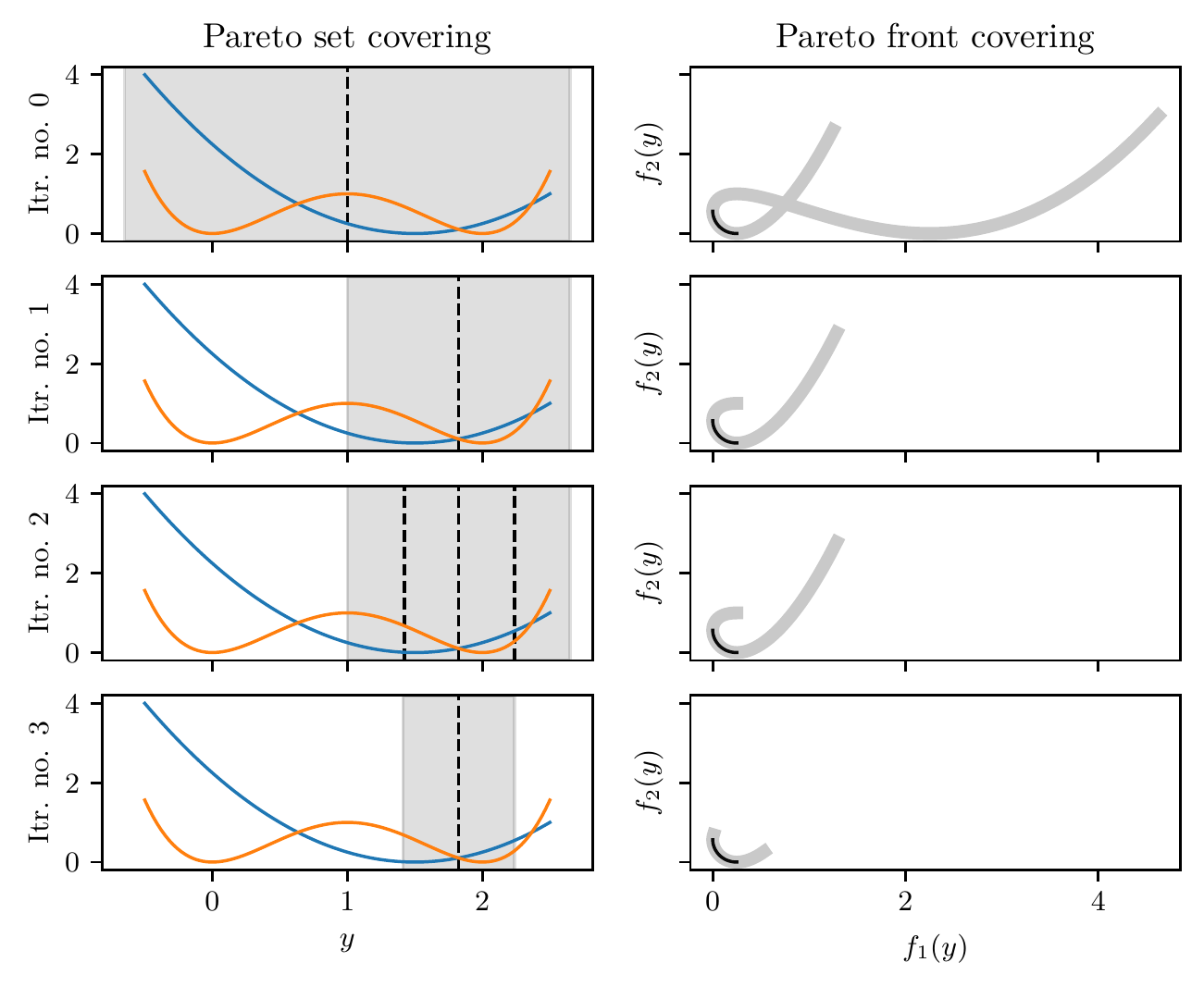}
    \caption{\new{First four iterations of the sampling algorithm demonstrated using Example~\ref{expl:pareto}. The dashed lines represent the boundaries of each box. In this example, each box is split in half before a non-dominance test is performed. The plots on the right show the images of the current box collections (gray/shaded) covering the true global Pareto front (black/solid).}}
    \label{fig:pareto-set-sampling-evolution}
\end{figure}

We will use a set-oriented method known as \emph{sampling algorithm}~\cite{Dellnitz2005}, which computes a box-covering of the Pareto set. The sampling algorithm is an iterative two-step process, \new{which in a first step subdivides a collection of boxes with respect to one coordinate and then discards every box that does not pass a set-wise non-dominance test, i.e., boxes that only contain dominated points. Figure~\ref{fig:pareto-set-sampling-evolution} shows how the boxes and the Pareto front covering evolve using the optimization problem introduced in Example~\ref{expl:pareto}. A complete description of the sampling algorithm can be found in Appendix~\ref{sec:appendix_sampling}. For further details, see also~\cite{Dellnitz2005, Schuetze2005, Schuetze2013, Peitz2017}.}

Each box can be efficiently represented by a center and a radius, so that all collections can be stored in a binary tree whose memory consumption grows linearly with~$n$, see~\cite{Dellnitz2001} for the MATLAB toolbox GAIO, which contains binary tree data structures and the algorithms for set-oriented calculations.

%%%%%%%%%%%%%%%%%%%%%%%%%%%%%%%%%%%%%%%%%%%%%%%%%%%%%%%%%%%%%%%%%%%%%%%%%%%

\section{Koopman-Based Surrogate Models with Control}
\label{sec:Koopman_Reduced_Models_with_Control}

For the following discussion, we consider non-deterministic dynamical systems with control of the form~\eqref{eq:SDE+control}. The aim of this section is to learn the Koopman generator associated with such systems. We consider two cases: (i) we assume that the control $u$ acts linearly on~\eqref{eq:SDE+control}, and (ii) we will consider the more general case, i.e., $u$ acts nonlinearly on~\eqref{eq:SDE+control}.

\subsection{Linear Case: Generator Interpolation}
\label{sec:Koopman_Generator_Interpolation}

We extend the recent result in~\cite{Peitz2020} for deterministic control-affine dynamical systems to \emph{non-deterministic} control-affine dynamical systems.

\begin{defn}
    A non-deterministic dynamical system of the from~\eqref{eq:SDE+control} is called \emph{control-affine} if
    \begin{align}
		b(X_t, u) &= b_0(X_t) + \sum_{i=1}^{d_u} b_i(X_t) \ts u_i, \label{eq:parameter_affine_drift} \\
		a(X_t, u) &= a_0(X_t) + \sum_{i=1}^{d_u} a_i(X_t) \ts u_i, \label{eq:parameter_affine_diffusion}
	\end{align}
    where $u_i \colon \mathbb{R}_{\ge 0} \to \mathbb{R}$ and $a = \sigma \ts \sigma^\top$.
\end{defn}
\begin{thm} \label{thm:controlaffine}
	Given a space of twice continuously differentiable functions and controls $u_1, u_2$, if the dynamics~\eqref{eq:SDE+control} are control-affine, then the Koopman generators are control-affine, i.e.,
    \begin{equation*}
        \mathcal{L}_{\alpha_1 u_1 + \alpha_2 u_2} = \mathcal{L}_0 + \alpha_1 \mathcal{A}_{u_1} + \alpha_2 \mathcal{A}_{u_2},
    \end{equation*}
    where $\mathcal{A}_u = \mathcal{L}_u - \mathcal{L}_0$ and $\alpha_1, \alpha_2 \in \mathbb{R}$.
\end{thm}

The proof of Theorem~\ref{thm:controlaffine} is analogous to the proof for deterministic control-affine systems, which can be found in, e.g.,~\cite{Peitz2020}. For the sake of completeness it is included in Appendix~\ref{sec:proof}.

\begin{expl} \label{expl:generator_interpolation}
Consider the nonlinear control-affine system
\begin{equation} \label{eq:control_affine_system}
    \ddt{t}
    \begin{bmatrix}
        x_1 \\
        x_2
    \end{bmatrix}
    = 
    \begin{bmatrix}
        (\gamma + g(u)) \ts x_1 \\
        \delta \ts (x_2 - x_1^2)
    \end{bmatrix},
\end{equation}
where \new{in this example} $u \equiv const$ and $g$ is a function. Defining $f_1 \coloneqq x_1$, $f_2 \coloneqq x_2$, and $f_3 \coloneqq x_1^2$, the system has a finite-dimensional, linear representation
\begin{equation} \label{eq:Koopman_inv_representation}
    \frac{\mathrm{d}}{\mathrm{d} t}
    \begin{bmatrix}
        f_1 \\
        f_2 \\
        f_3
    \end{bmatrix}
    =
    \underbracedmatrix{
        \gamma + g(u) & 0 & 0 \\
        0 & \delta & -\delta \\
        0 & 0 & 2 \ts (\gamma + g(u))
    }{=L_u}
    \begin{bmatrix}
        f_1 \\
        f_2 \\
        f_3
    \end{bmatrix},
\end{equation}
where $L_u$ is a finite-dimensional representation of the generator $\mathcal{L}_u$. Following Theorem~\ref{thm:controlaffine}, \new{if $g$ is linear, $L_u$} can be split into 
\begin{equation*}
    L_u = 
    \underbracedmatrix{
        \gamma & 0 & 0 \\
        0 & \delta & -\delta \\
        0 & 0 & 2 \ts \gamma
    }{=L_0}
    +
    \underbracedmatrix{
        g(u) & 0 & 0 \\
        0 & 0 & 0 \\
        0 & 0 & 2 \ts g(u)
    }{=A_u},
\end{equation*}
where $A_u$ is a finite-dimensional representation of $\mathcal{A}_u$ and linear with respect to $u$. Thus, $L_u$ is control-affine. \exampleSymbol
\end{expl}

\subsection{Nonlinear Case: State Augmentation}

To motivate what follows, take again a look at Example~\ref{expl:generator_interpolation} and assume that $g$ is nonlinear. The system~\eqref{eq:control_affine_system} can still be represented as in~\eqref{eq:Koopman_inv_representation}, however, nonlinearity of $g$ prohibits to use Theorem~\ref{thm:controlaffine} as the operator $A_u$ is not linear with respect to $u$. A workaround is to augment the system so that the control $u$ is represented as an additional state, i.e., the augmented system is given by
\begin{equation*}
    \ddt{t}
    \begin{bmatrix}
        x_1 \\
        x_2 \\
        u
    \end{bmatrix}
    = 
    \begin{bmatrix}
        (\gamma + g(u)) \ts x_1 \\
        \delta \ts (x_2 - x_1^2) \\
        0
    \end{bmatrix}.
\end{equation*}
A closed representation of the augmented system as in~\eqref{eq:Koopman_inv_representation} for the non-augmented system~\eqref{eq:control_affine_system}, however, is not possible any longer. We will now show how gEDMD can be extended to control inputs. The idea \new{of augmenting the state space} has also been used in other data-driven methods such as \new{DMDc~\cite{Proctor2016} or EDMDc~\cite{Korda2018}, both of which provide approximations of the Koopman operator using linear and nonlinear basis functions, respectively, or SINDYc~\cite{Brunton2016a}, which can be used to find the governing equations of a dynamical system with control.}

For a general system~\eqref{eq:SDE+control}, we assume that we have a set of $ m $ measurements of the augmented system state $ \{\ts [x_l, u_l] \ts\}_{l=1}^m $ as well as $ \{\ts \dot u_l \ts\}_{l=1}^m $, which is trivial as we have full access to $u$, i.e., we can choose its values and thus know its derivative. Let $ [x_l, u_l] = [x_{l_1}, \dots, x_{l_d}, u_{l_1}, \dots, u_{l_{d_u}}] \in \mathbb{R}^{d+d_u}$. To simplify the notation, let $\bar x_l \coloneqq [x_l, u_l]$. Additionally, assume that we have a set of the augmented drift $ \{\ts [b(\bar x_l), \dot u_l] \ts\}_{l=1}^m $ as well as the augmented diffusion $ \{\ts \sigma(\bar x_l, u_l) \ts\}_{l=1}^m $, which is defined as
\begin{equation*}
    \sigma(\bar x_l, u_l) \coloneqq
    \begin{bmatrix}
        \sigma(\bar x_l) & \sigma_{12}(u_l) \\
        \sigma_{21}(u_l) & \sigma_{22}(u_l)
    \end{bmatrix} \in \mathbb{R}^{(d+d_u) \times (d+d_u)},
\end{equation*}
where $\sigma_{12}(u_l) \in \mathbb{R}^{d \times d_u}$, $\sigma_{21}(u_l) \in \mathbb{R}^{d_u \times d}$ and $\sigma_{22}(u_l) \in \mathbb{R}^{d_u \times d_u}$ denote the diffusion terms of the control $u$ (typically zero). We assume that either both drift $b(\bar x_l)$ and diffusion $\sigma(\bar x_l)$ are known or that they can be estimated pointwise. Let $ \{\ts \psi_k \ts\}_{k=1}^\ell $ denote the set of basis functions and define
\begin{equation*}
	\mathrm{d}\psi_k(\bar x) \coloneqq (\mathcal{L} \psi_k)(\bar x) = \sum_{i=1}^{d+d_u} b_i(\bar x) \ts \pd{\psi_k}{\bar x_i}(\bar x) + \frac{1}{2} \sum_{i=1}^{d+d_u} \sum_{j=1}^{d+d_u} a_{ij}(\bar x) \ts \pd{^2 \psi_k}{\bar x_i \ts \partial \bar x_j}(\bar x),
\end{equation*}
where $a = \sigma(\bar x_l, u_l) \ts \sigma(\bar x_l, u_l)^\top$. The partial derivatives of the basis functions can be precomputed analytically. We can then compute the gEDMD matrices $\Psi_X$ and $\mathrm{d}\Psi_X$ in~\eqref{eq:gEDMDmatrices} for $\bar x_1, \dots, \bar x_m$ and solve the minimization problem $\min \Vert \mathrm{d}\Psi_X - M \Psi_X \Vert_F$ to obtain a finite-dimensional empirical estimate of the generator $L = M^\top$. 

Note that we make three important restrictions on the feasible controls: (i)~the controls are given by families of solutions of ODEs, (ii)~the controls act on long time scales only (e.g., bang-bang controls are not allowed), and (iii)~the controls are representable in terms of the basis functions. However, these constraints do not limit our approach for the following reasons: (i)~in most scenarios, controls are typically given in such a form, (ii)~rapidly changing controls are unrealistic in the context of ABMs, and (iii)~insufficient and inappropriately chosen basis functions lead to inaccurate results.

\begin{figure}
    \centering
    \begin{subfigure}[t]{0.49\textwidth}
        \includegraphics[width=\textwidth]{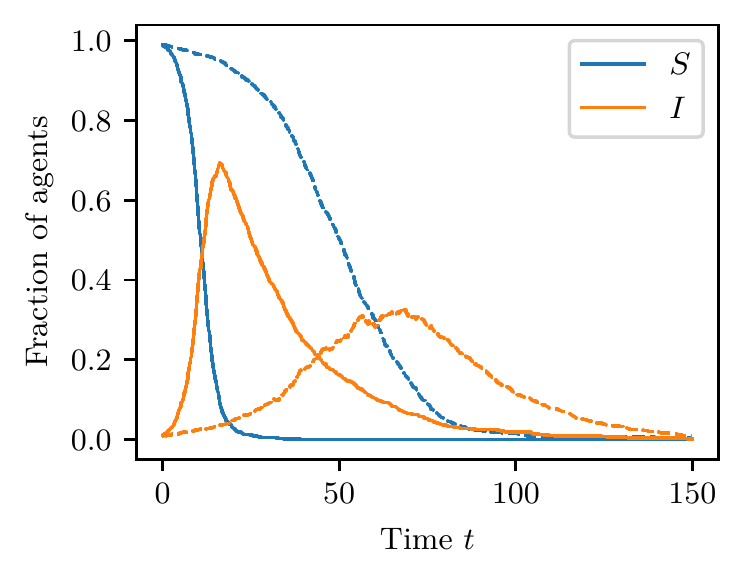}
        \caption{}
    \end{subfigure}
    \begin{subfigure}[t]{0.49\textwidth}
        \includegraphics[width=\textwidth]{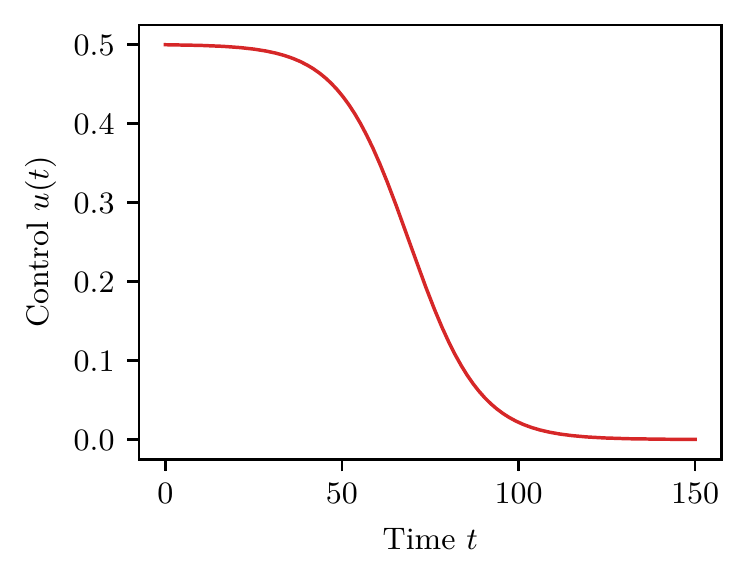}
        \caption{}
    \end{subfigure}
    \caption{(a)~Trajectories of susceptible and infected agents as a fraction of the population for the SIR model with (dashed) and without (solid) control $u$ and (b)~control $u(t)$ applied to SIR model.}
    \label{fig:SIR_input_example}
\end{figure}

\begin{expl} \label{expl:SIR-ABM}
    To illustrate \new{how the state augmentation works}, consider a stochastic SIR model to simulate an infectious disease and its containment with non-phar\-ma\-ceu\-ti\-cal interventions. For simplicity, we consider a population of size $N$ and assume that an infection occurs at rate $\beta(u)$, which depends on the disease itself and the containment strategy $u(t)$. Once infected, individuals recover at rate $\gamma$ and cannot be re-in\-fec\-ted. Defining the system state $X_t \coloneqq [S_t, I_t]^\top$, where $S_t$ and $I_t$ denote the numbers of susceptible and infected individuals, respectively, the stochastic SIR model is given by
    \begin{align*}
        \mathrm{d}S_t &= - \beta(u) \ts \frac{S_t I_t}{N} \ts \mathrm{d}t - \frac{1}{N} \sqrt{\beta(u) \ts S_t I_t} \ts \mathrm{d}W_1(t), \\
        \mathrm{d}I_t &= \left[\beta(u) \ts \frac{S_t I_t}{N} -\gamma \ts I_t\right] \mathrm{d}t + \frac{1}{N} \sqrt{\beta(u) \ts S_t I_t} \ts \mathrm{d}W_2(t) - \sqrt{\gamma \ts N^{-1} \ts I_t} \ts \mathrm{d}W_3(t).
    \end{align*}
    The number of recovered individuals is given by $R_t = N - S_t - I_t$ since the size of the population is assumed to be constant.
    
    It has been known -- not only since the COVID-19 pandemic -- that non-phar\-ma\-ceu\-ti\-cal interventions such as social distancing, wearing face masks, or lockdowns can slow or prevent the spread of infectious diseases. Such interventions can be expected to have a nonlinear effect on the dynamics because the probability of infection changes nonlinearly with the viral load in the air to which a person is exposed~\cite{Cheng2021}. For this toy model we assume that $u$ acts quadratically via $\beta(u) = \tilde\beta (1 - u(t))^2$ for time $t \ge 0$.  The system is not control-affine. Thus, Theorem~\ref{thm:controlaffine} cannot be applied. State augmentation, on the other hand, can be applied. We assume that the control $u$ is a continuously relaxed intervention of the form $u(t) = A \ts Q/(1 + e^{B \ts t})$. We set $\tilde\beta = 0.5$, $\gamma = 0.05$, $A=0.5$, $Q=1\ts 000$, and $B=0.1$. Figure~\ref{fig:SIR_input_example}(a) shows the trajectories for susceptible and infected individuals as a fraction of the population without (solid) and with control $u$ (dashed). The control $u(t)$ is shown in Figure~\ref{fig:SIR_input_example}(b). 
    
    To compute a representation of the Koopman generator, we augment the system state by the control $u$ and generate $m = 1\ts 000$ uniformly distributed training points $\{\ts [x_l, u_l] \ts\}$, where $x_l \coloneqq [S_t/N, I_t/N]_l^\top$. We use exact values for the augmented drift $ \{\ts [b(x_l, u_l), \dot u_l] \ts\}_{l=1}^m $ as well as the augmented diffusion $ \{\ts \sigma(x_l, u_l) \ts\}_{l=1}^m$. For the basis functions we choose monomials up to degree 5. The system is correctly identified with drift and diffusion terms given by
    \begin{align*}
        b(x, u) &=
        \begin{bmatrix}
            -\tilde\beta \ts x_1 x_2 (1-u)^2 \\
            \tilde\beta \ts x_1 x_2 (1-u)^2 -\gamma \ts x_2 \\
            B \ts u + \frac{B}{A} u^2
        \end{bmatrix} \\
        \intertext{and} a(x, u) &= \frac{1}{N}
        \begin{bmatrix}
            \tilde\beta \ts x_1 x_2 (1-u)^2 & -\tilde\beta \ts x_1 x_2 (1-u)^2 & 0\\
            -\tilde\beta \ts x_1 x_2 (1-u)^2 & \tilde\beta \ts x_1 x_2 (1-u)^2 + \gamma \ts x_2 & 0 \\
            0 & 0 & 0
        \end{bmatrix}.
    \end{align*}
    Since in general we are only interested in identifying the dynamics of $X_t = [S_t, I_t]^\top$ and not the dynamics of the control $u$ (which is actually known), we could relax the third restriction concerning the representation in terms of the basis function and neglect the columns in $L$ corresponding to $u$. For further details on system identification and model reduction using gEDMD, see~\cite{Klus2020}. \exampleSymbol
\end{expl}

%%%%%%%%%%%%%%%%%%%%%%%%%%%%%%%%%%%%%%%%%%%%%%%%%%%%%%%%%%%%%%%%%%%%%%%%%%%

\section{A Multi-Objective Optimization Approach for ABMs}
\label{sec:MOOP4ABMS}

The aim in this section is to find the Pareto set of $n$ controls $u_i$ that optimize $k$ objective functions $f_i \colon \mathcal{R} \to \mathbb{R}$ in \eqref{eq:MOOP} depending on the state of an ABM. We assume that each function $f_i$ can be written as
\new{
\begin{equation*}
	f_i(X, u) = \mathbb{E}\left[ \int_{t_0}^{t_1} r_i (X(t; u), u(t)) \, \mathrm{d}t + s_i(X(t_1, (u(t_1)))) \right], \quad t_0 < t_1,
\end{equation*}}
where $r_i$ and $s_i$ are \new{the running and terminal cost} functions depending on state \new{$X(t; u)$} of an ABM and  $u$ is the control. \new{An optimal control $u^\star$ is given by a solution of the multi-objective optimization problem
\begin{equation*}
    \underset{u \in U}{\text{min}} ~ F(X, u),
\end{equation*}
where $F(X, u) = \left[f_1(X, u), \dots, f_k(X, u)\right]^\top$ and $U = \{u \colon \mathbb{R}_{\ge 0} \to \mathbb{R}^{d_u} \mid u(\cdot) \ \text{measureable}\}$ the space of admissible controls.} Three major issues occur when optimizing ABMs:
\begin{enumerate}[label=(\roman*), itemsep=0ex, topsep=0.5ex]
	\item \textbf{Computational complexity:} Most ABMs involve thousands of agents and simulate complex interactions between agents, often at high time resolution, making the calculation of single trajectories computationally expensive.
    \item \textbf{Stochasticity:} Most ABMs are inherently stochastic so that many independent simulations are required to compute statistical properties such as $ \mathbb{E}[\,\cdot\,] $.
    \item \textbf{Discontinuity:} Most ABMs are inherently discontinuous since agents make discrete decisions, and thus single trajectories are not differentiable with respect to $u$.
\end{enumerate}
By replacing the ABM with a surrogate model that is less expensive to evaluate and directly approximates expectations, we address the problems of computational complexity and stochasticity. It has recently been shown in~\cite{Niemann2022} that under certain conditions it is possible to obtain reduced models, \new{given by differential equations,} from ABM data that accurately approximate the \new{aggregated} dynamics of ABMs, i.e., the collective behavior of the agents. \new{That is, instead of dealing with the ABM state space whose size grows combinatorially as $d^N$, where $N$ denotes the number of agents and $d$ the number of states each agent can have, the space of the surrogate model is at most $d$-dimensional (with $d \ll N$).} Together with the results presented in Section~\ref{sec:Koopman_Reduced_Models_with_Control}, we now construct computationally cheap surrogate models that preserve features relevant for optimizing ABMs.

The sampling algorithm introduced in Section~\ref{sec:MOOP} avoids the problem of discontinuity and the resulting nonexistence of derivatives since it is derivative-free, i.e., it relies only on function evaluations. We will demonstrate in Sections~\ref{sec:VM_MOO} and~\ref{sec:GERDA_MOO} how data-driven reduced models of ABMs can be effectively used as surrogate models in multi-objective optimization problems. 

\new{
\begin{rem}
	For the sake of illustration, we will only consider constant controls. Nevertheless, the approach can also be applied to time-varying controls such as different levels of curfews that change, e.g., every week. Additionally, having a reduced dynamical model based on the Koopman formalism opens up the possibility to use existing methods for controlling linear problems.
\end{rem}}

\new{We will now introduce two ABMs which will be used as examples to demonstrate the different approaches described in Section \ref{sec:Koopman_Reduced_Models_with_Control}.}

\begin{figure}
    \centering
    \begin{subfigure}[t]{0.49\textwidth}
        \includegraphics[width=\textwidth]{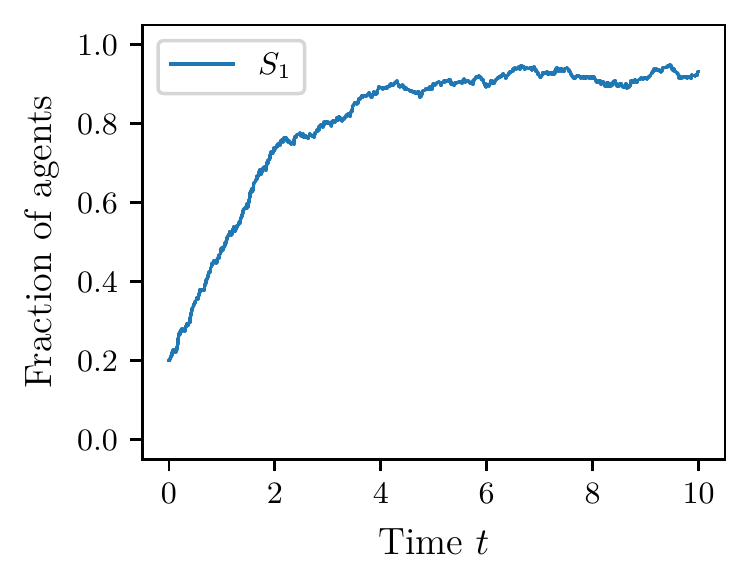}
        \caption{}
    \end{subfigure}
    \begin{subfigure}[t]{0.49\textwidth}
        \includegraphics[width=\textwidth]{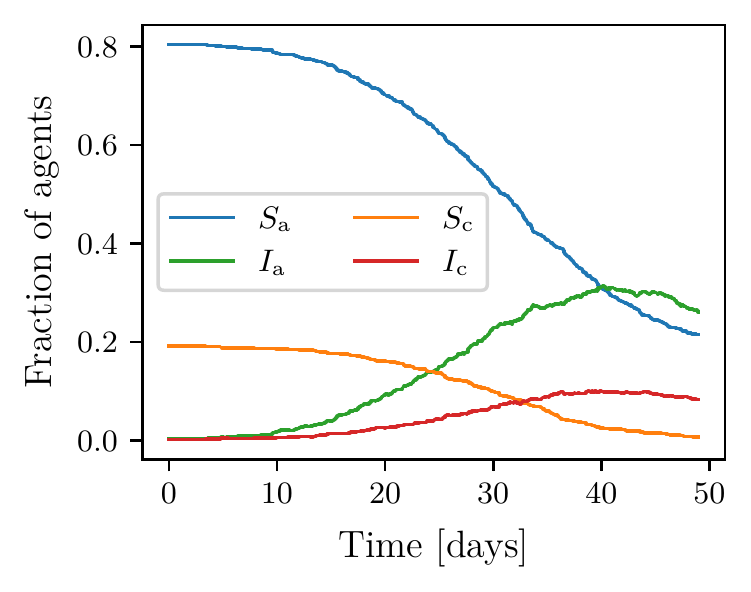}
        \caption{}
    \end{subfigure}
    \caption{Aggregated trajectories of (a)~the voter model and (b)~the GERDA model for the parameters given in Table~\ref{tab:parameters}.}
    \label{fig:GERDA_and_VM}
\end{figure}

\paragraph{\new{The Voter Model.}} The voter model was first introduced in~\cite{Holley1975} and has since not only been used in the social context but also in many other contexts, such as chemical systems or ant colonies~\cite{HerreriasAzcue2019, Biancalani2014, Ohkubo2008}. In its basic definition, $N$ identical agents interact with each other at any time in a given network and change their opinions based on (stochastic) transition rules. In this work, we consider the model with two opinions and two transition rules. Given two agents with opinions $S_i \neq S_j$, the first transition rule is given by
\begin{equation*}
	R_{ij} \colon ~~ S_i + S_j \xmapsto{\gamma_{ij}} 2 \ts S_j,
\end{equation*}
which is a second-order transition, meaning that an agent with opinion $S_i$ adopts another agent's opinion $S_j$ at rate $\gamma_{ij}$. The second transition rule is a first-order transition of the form
\begin{equation*}
	R_{ij}' \colon ~~ S_i \xmapsto{\gamma_{ij}'} S_j,
\end{equation*}
where an agent changes its opinion independently of all other agents at rate $\gamma'_{ij}$. Gillespie's stochastic simulation algorithm~\cite{Gillespie1976}, which constructs exact realizations in con\-tin\-u\-ous time, can be used to simulate the ABM. We choose $N=500$ agents and set $\gamma_{12} = 1$, $\gamma_{21} = 2$, and $\gamma_{12}' = \gamma_{21}' = 0.1$. A single aggregated simulation is shown in Figure~\ref{fig:GERDA_and_VM}(a). That is, we visualize the number of agents with opinion $S_1$ at time $t$. The trajectory for $S_2$ is omitted due to conservation of $N$, i.e., $S_2 = N - S_1$. For further details, especially with respect to convergence to ordinary or stochastic differential equations, we refer the reader to~\cite{Kurtz1978}.

\paragraph{\new{The Georeferenced Demographic Agent-Based Model.}} The second example is the individualized \textbf{GE}o\-\textbf{R}ef\-er\-enced \textbf{D}e\-mo\-graph\-ic \textbf{A}gent-based mod\-el (GERDA) for the transmission of SARS-CoV-2 and the disease dynamics of COVID-19~\cite{Goldenbogen2022}. Using detailed location data, a demographically matched population, and realistic, hourly schedules for each agent, this model simulates contacts between people in the given locations and the resulting infection events. Data from multiple towns, e.g., Tepoztlán (Mexico), Zikhron Ya'akov (Israel), or Gangelt (Germany) have been used to calibrate the GERDA model and analyze the respective infection dynamics. In this paper, the calibrated data of the German municipality of Gangelt is used. For more information, see the authors' original publication. The model is implemented in Python and is available at \url{https://ford.biologie.hu-berlin.de/jwodke/corona_model}.

To simplify the disease dynamics, which originally distinguishes between agents that are susceptible, diagnosed, hospitalized, in intensive care, and recovered or deceased, we summarize and consider only the classical three compartments $S$ (susceptible), $I$ (infected), and $R$ (recovered) for the aggregated dynamics. Since in the subsequent analysis we consider homeschooling and \new{working from home} as controls, and since a separation by age is also justified by different infection rates (see, e.g.,~\cite{She2020}), we group the agents by children (age 0 to 18 years) and adults (19\texttt{+}) and use the subscripts $\mathrm{c}$ and $\mathrm{a}$, respectively. Figure~\ref{fig:GERDA_and_VM}(b) shows an aggregated trajectory of a single simulation of GERDA in terms of numbers of susceptible, infected, and recovered adults and children for the parameters given in Table~\ref{tab:parameters}. The trajectories for $R_\mathrm{a}$ and $R_\mathrm{c}$ are omitted due to the conservation of the number of adults $N_\mathrm{a}$ and children $N_\mathrm{c}$.

\subsection{Multi-Objective Optimization of the Voter Model}
\label{sec:VM_MOO}

\begin{figure}
	\centering
    \begin{subfigure}[t]{0.49\textwidth}
		\centering
		\includegraphics[width=\textwidth]{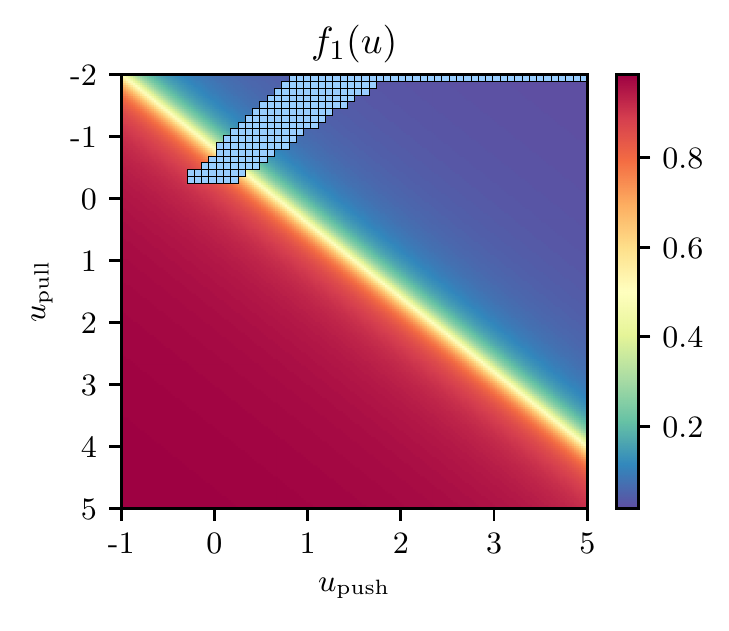}
		\caption{}
	\end{subfigure}
    \begin{subfigure}[t]{0.49\textwidth}
		\centering
		\includegraphics[width=\textwidth]{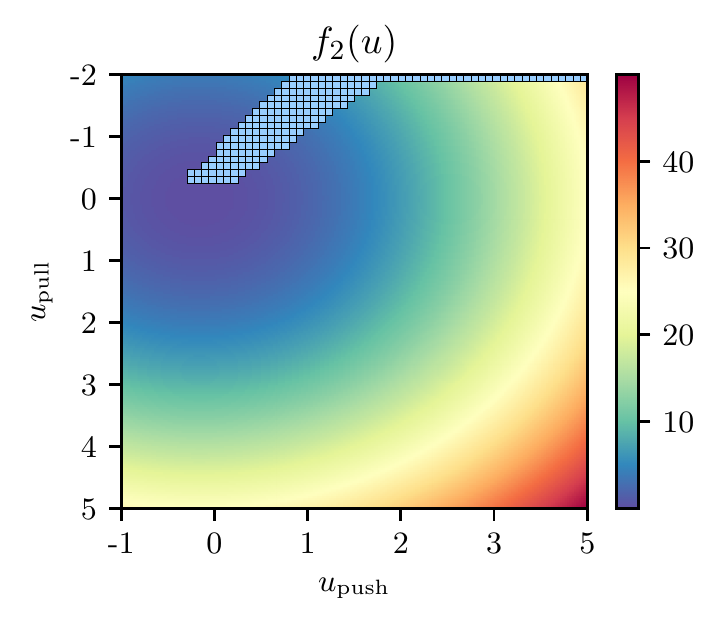}
		\caption{}
	\end{subfigure}
	\caption{Pareto set covering (light blue boxes) after 12 iterations plotted on (a)~objective function~\eqref{eq:VM_expensive} and (b)~objective function~\eqref{eq:VM_cheap} at time $t=10$ depending on the control $ [u_\text{push}, u_\text{pull}]^\top \in \mathcal{R}$. The yellow line in~(a) marks the area of interest where the majority switches from $S_1$ to $S_2$.}
	\label{fig:VM_objective_functions}
\end{figure}

In order to win over as many voters as possible to opinion $S_2$, we want to find combinations of push and pull arguments that simultaneously move voters away from opinion $S_1$ and towards opinion $S_2$. We represent this behavior in terms of the control $u = [u_\text{push}, u_\text{pull}]^\top \in \mathcal{R}$, where $u_\text{push}$ acts on the transition rate $\gamma_{12}$ from $S_1$ to $S_2$ and $u_\text{pull}$ acts analogously on $\gamma_{21}$. We assume that $u$ acts linearly on the dynamics given by the ABM, i.e., $\gamma_{12}(u) = \tilde\gamma_{12} + u_\text{push}$ and $\gamma_{21}(u) = \tilde\gamma_{21} + u_\text{pull}$, and set up a multi-objective optimization problem of the form~\eqref{eq:MOOP} with feasible decision space $\mathcal{R} \coloneqq [-1, 5] \times [-2, 5]$ and objective function $F \colon \mathcal{R} \to \mathbb{R}^2$ with
\begin{align}
    f_1(u) &= \mathbb{E}[X_1(t; u) / N], \label{eq:VM_expensive}\\
    f_2(u) &= u_\text{push}^2 + u_\text{pull}^2, \label{eq:VM_cheap}
\end{align}
where $X_1(t; u)$ denotes the number of agents believing in opinion $S_1$ at time $t$ depending on control $u \in \mathcal{R}$. The second objective is chosen such that any effort to influence the agent's opinion is costly. We evaluate the objective~\eqref{eq:VM_expensive} at time $t=10$. Figures~\ref{fig:VM_objective_functions}(a) and~(b) show the objective functions~\eqref{eq:VM_expensive} and~\eqref{eq:VM_cheap}, respectively.

Even though this ABM is computationally inexpensive in comparison to other ABMs, it is still inefficient to optimize it directly. Therefore, we first construct a surrogate model using gEDMD, and then solve the multi-objective optimization problem. As the control $u$ acts linearly on the dynamics, Theorem~\ref{thm:controlaffine} guarantees that we can obtain a surrogate model for varying control $u$ using interpolated Koopman generators. In fact, for sufficiently large $N$ and due to conservation of $N$, the drift and diffusion terms identified by gEDMD correspond to the well-known SDE approximation by Kurtz~\cite{Kurtz1978} (also known as the chemical Langevin equation~\cite{Gillespie2000}) satisfying
\begin{align*}
		\mathrm{d}C_t  &=  \Big[ ((\gamma_{21} + u_\text{pull}) - (\gamma_{12} + u_\text{push})) \ts C_t \ts (1-C_t)-\gamma'_{12} \ts C_t + \gamma'_{21}(1-C_t) \Big] \mathrm{d}t \\
		&+ \frac{1}{\sqrt{N}} \Big[ -\sqrt{(\gamma_{12} + u_\text{push}) \ts C_t \ts (1-C_t)} \ts \mathrm{d}W_1(t) +\sqrt{(\gamma_{21} + u_\text{pull}) \ts C_t \ts (1-C_t)}\ts  \mathrm{d}W_2(t) \Big] \\
		&+ \frac{1}{\sqrt{N}} \Big[ -\sqrt{\gamma'_{12} \ts C_t} \ts  \mathrm{d}W_3(t) + \sqrt{\gamma'_{21}(1-C_t)} \ts  \mathrm{d}W_4(t) \Big],
	\end{align*}
where $C_0 = \lim_{N \to \infty} X_1(0; u)/N$. Drift and diffusion terms can be written as~\eqref{eq:parameter_affine_drift} and~\eqref{eq:parameter_affine_diffusion}. Thus, Theorem~\ref{thm:controlaffine} holds and we learn four generator approximations to cover the feasible decision space $\mathcal{R}$. We use the four controls marking the vertices of $\mathcal{R}$ to construct the surrogate model. The matrix representation of the Koopman generator approximation for the surrogate is then given by $L_u = \sum_{i = 1}^{4} \alpha_i \ts L(u_i)$, where $\sum_{i = 1}^{4} \alpha_i = 1$ with $\alpha_i \in [0, 1]$ and $L(u_1) = L([0, 0])$, $L(u_2) = L([1, 0])$, $L(u_3) = L([0, 1])$ and $L(u_4) = L([1, 1])$ denote the matrix representations of the Koopman generator approximation for the dynamics with control $u_i$. To train the reduced models at the vertices, i.e., learn the matrix representations $L(u_i)$, we sample the ABM at 100 uniformly distributed points, each having 100 Monte Carlo simulations to calculate pointwise drift and diffusion estimates using the Kramers--Moyal formulae. For details on Kramers--Moyal expansions, see~\cite{Risken1996}.

\begin{figure}
	\centering
	\begin{subfigure}[t]{0.49\textwidth}
		\centering
		\includegraphics[width=\textwidth]{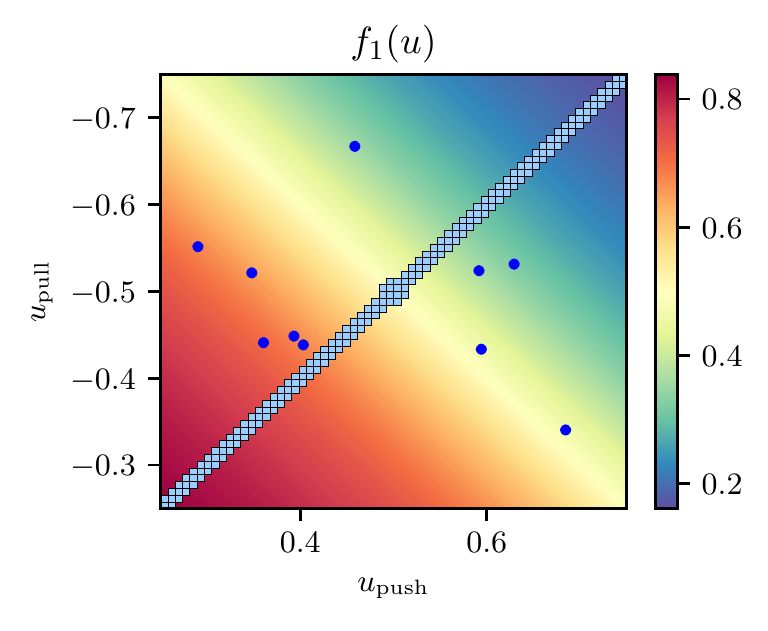}
		\caption{}
	\end{subfigure}
	\hfil
	\begin{subfigure}[t]{0.49\textwidth}
		\centering
		\includegraphics[width=\textwidth]{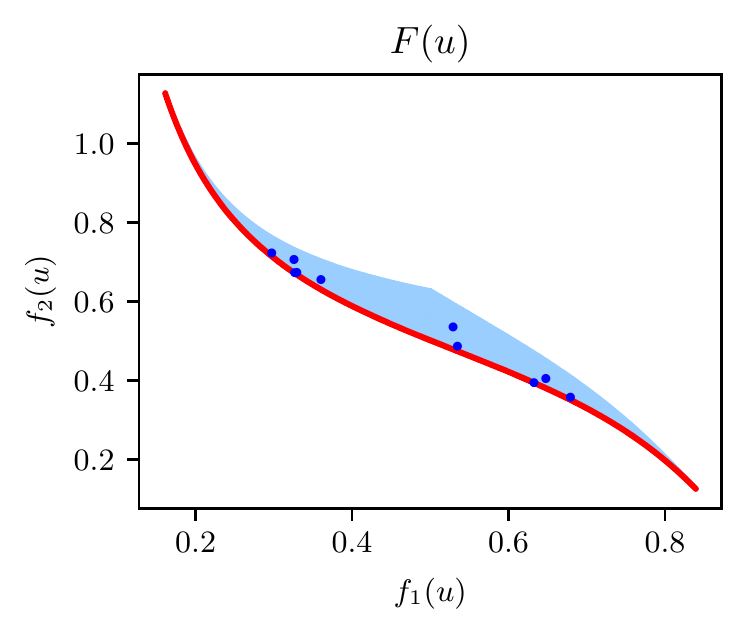}
		\caption{}
	\end{subfigure}
	\begin{subfigure}[t]{\textwidth}
		\centering
		\includegraphics[width=\textwidth]{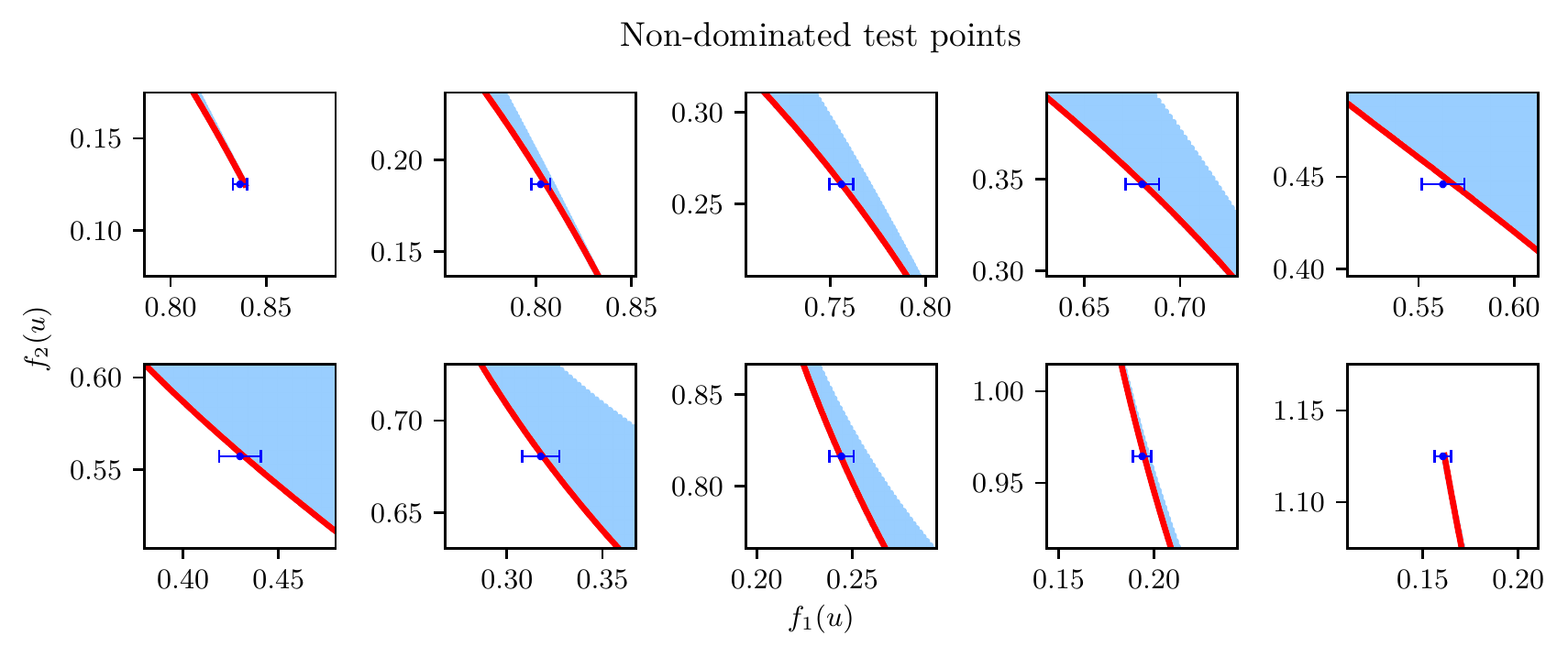}
        \includegraphics[width=\textwidth]{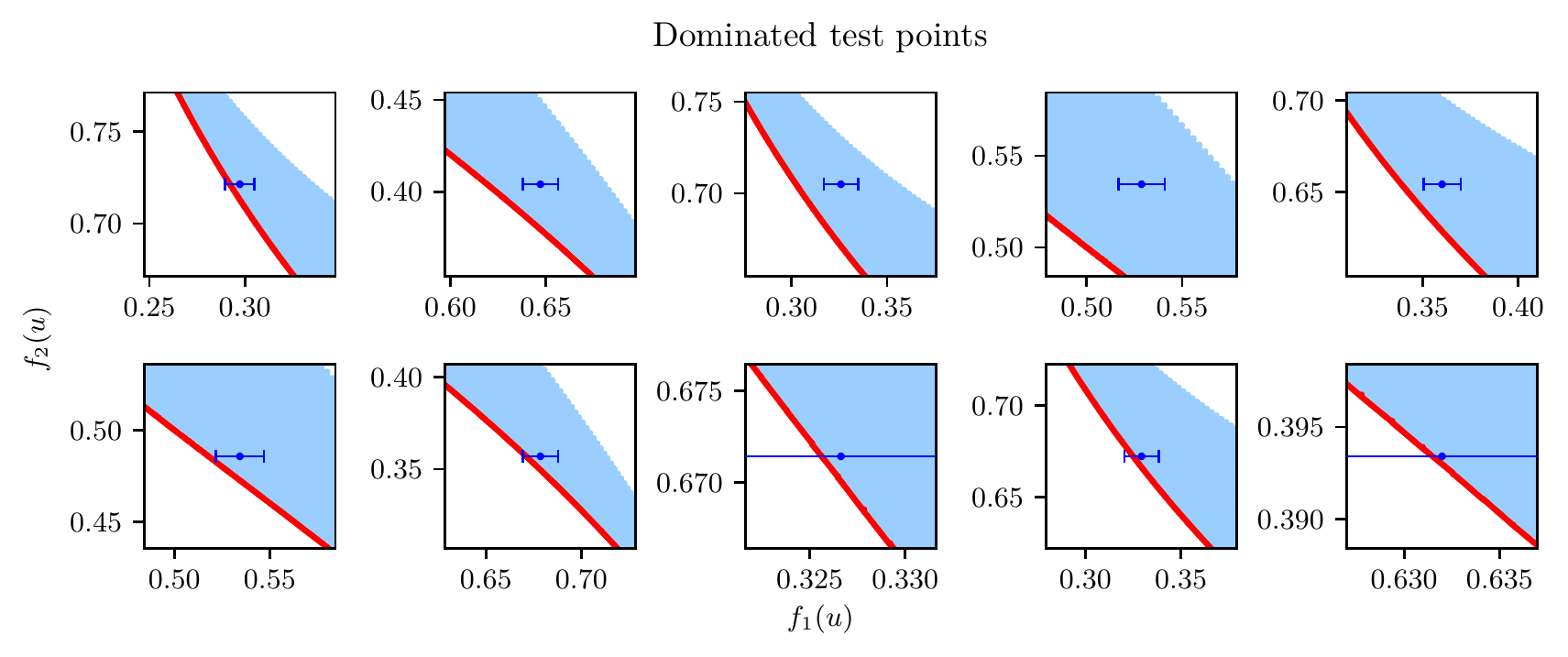}
		\caption{}
	\end{subfigure}
	\caption{(a)~Pareto set covering (light blue boxes) and test points (blue/dots) for refined feasible decision space $\mathcal{R}^\star$. (b)~Pareto front (red/solid), image of $\mathcal{R}^\star$ (light blue area) and test points under the objective function $F$. (c)~Close-ups of randomly chosen test points, evaluated via the ABM. While the non-dominated test points are covered by the boxes in (a) and thus are Pareto-optimal points, the dominated test points are not. Error bars indicate a 99.9~\% confidence level.}
	\label{fig:pareto_front_refined}
\end{figure}

Figure~\ref{fig:VM_objective_functions} shows the Pareto set covering obtained after 12 iterations for the feasible decision space $\mathcal{R}$ together with the objectives~\eqref{eq:VM_expensive} and~\eqref{eq:VM_cheap}. In order to achieve an absolute majority for opinion $S_2$, we refine the decision space to $\mathcal{R}^\star = [0.25, 0.75] \times [-0.75, -0.25]$. Figure~\ref{fig:pareto_front_refined} shows in~(a) the computed covering of the Pareto set for $\mathcal{R}^\star$ and in (b)~the image of $\mathcal{R}^\star$ under objective function $F$ as a light blue area as well as the approximated Pareto front (red/solid)\new{, which is obtained by mapping the center of each box covering the Pareto set in~(a)}. To verify that the surrogate model approximates the dynamics of the ABM sufficiently well and consequently is also capable of approximating the Pareto set, we randomly choose some test points (blue/dots) for which we evaluate the objective~\eqref{eq:VM_expensive} from 100 Monte Carlo simulations using the ABM. Figure~\ref{fig:pareto_front_refined}(c) shows close-ups of these test points. We observe that the test points that are covered by the boxes in Figure~\ref{fig:pareto_front_refined}(a) (and are therefore not visible) are non-dominated points. These points are mapped to the Pareto front computed via the surrogate model. The visible (blue marked) test points (in general, points that are not covered by boxes in Figure~\ref{fig:pareto_front_refined}(a)) are on the right side of the Pareto front and thus dominated. Thus, the box covering computed via the data-driven surrogate model actually approximates the Pareto-optimal set of the ABM as it captures the dynamical behavior of the ABM under control $u$ sufficiently well. The error bars indicate a confidence level of 99.9~\%. 

\subsection{Multi-Objective Optimization of the \new{GERDA} Model}
\label{sec:GERDA_MOO}

The design of optimal containment strategies for a complex model like GERDA has recently been discussed in~\cite{Niemann2023}. Although the considered objective functions are only prototypical, they show that the design of optimal containment strategies using ABMs crucially depends on the level of detail of the model. However, even though the main objective in~\cite{Niemann2023} consists of three sub-objectives (i.e., health, society, and economy), multi-objective optimization is not considered. We now show how a multi-objective optimal containment strategy can be found using a suitable reduced model of GERDA.

Following~\cite{Niemann2023}, we consider \emph{constant} controls $u = [u_\mathrm{s}, u_\mathrm{w}]^\top \in \mathcal{R}$ and set the feasible decision space for our multi-objective optimization problem~\eqref{eq:MOOP} to $\mathcal{R} = [0, 1] \times [0, 0.8]$ and, for better comparability, choose the same objectives and parameters as in~\cite{Niemann2023}, i.e., $F \colon \mathcal{R} \to \mathbb{R}^2$ with
\begin{align} 
        f_1(u) &= \mathbb{E} \left[ \int_0^T (I(t; u)/N + \exp(10(I(t; u) - I_\mathrm{max})/N)) \, \mathrm{d}t \right], \label{eq:GERDA_health_objective} \\
        f_2(u) &= \int_0^T  \left[ u_\mathrm{s}(t)^2 -\log(u_\mathrm{w}^\mathrm{max} - u_\mathrm{w}(t)) \right]  \mathrm{d}t, \label{eq:GERDA_economic_objective}
\end{align}
where $f_1$ and $f_2$ denote the objectives related to health and society/economy, respectively. Further, let $I(t; u)$ denote the number of infected agents, $I_\mathrm{max}$ and $u_\mathrm{w}^\mathrm{max}$ threshold values, $u_\mathrm{s}$ and $u_\mathrm{w}$ the controls representing the fraction of school and work closures (where $u_\mathrm{s} = 0$ is interpreted as no schools being closed and $u_\mathrm{s} = 1$ as all being closed; analogously for $u_\mathrm{w}$). Note that we do not consider the objectives for society and economy individually as this would result in a trivial multi-objective optimization problem. In fact, the whole domain $\mathcal{R}$ would be Pareto-optimal since in this case each objective increases in a different direction. 

\begin{rem}
    Equation~\eqref{eq:GERDA_health_objective} accounts for the direct negative consequences of infections, which is assumed to exhibit an approximate linear relationship with the number of infected individuals, as well as for the social impact, which increases significantly when the capacity of the health care system to treat severely ill individuals is exceeded. This is quantified by a threshold of infected individuals, denoted as $I_\mathrm{max}$. Equation~\eqref{eq:GERDA_economic_objective} represents the direct costs of both controls $u_\mathrm{s}$ and $u_\mathrm{w}$. The economic costs become $+\infty$ when the \new{homeworking} rate approaches an upper bound $u_\mathrm{w}^\mathrm{max} < 1$ in case that too many workplaces are closed or employees are working from home.
\end{rem}

\begin{figure}
    \centering
    \includegraphics[width=\textwidth]{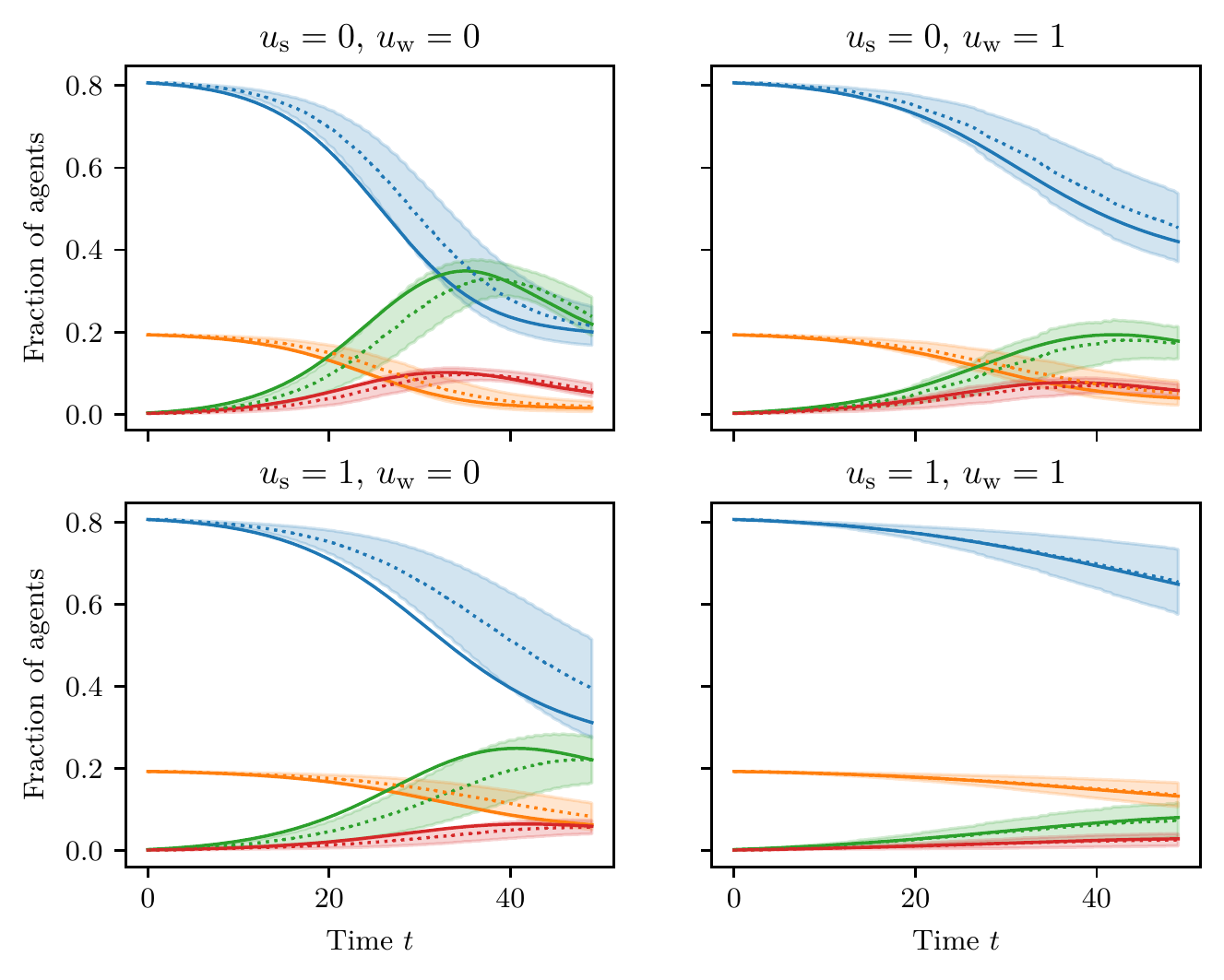}
    \caption{Mean aggregated trajectories for GERDA (dotted) and the reduced model (solid) of the fraction of susceptible (blue for adults, orange for children) and infected agents (green for adults, red for children) under different controls $u=[u_\text{s}, u_\text{w}]^\top$. The shaded areas indicate the standard deviations of the GERDA model. The ABM data is estimated from 1\ts 000 independent simulations using a fixed world and the parameters in Table~\ref{tab:parameters}.}
    \label{fig:GERDA_reduced}
\end{figure}

\begin{figure}
    \centering
    \begin{subfigure}[t]{0.49\textwidth}
        \includegraphics[width=\textwidth]{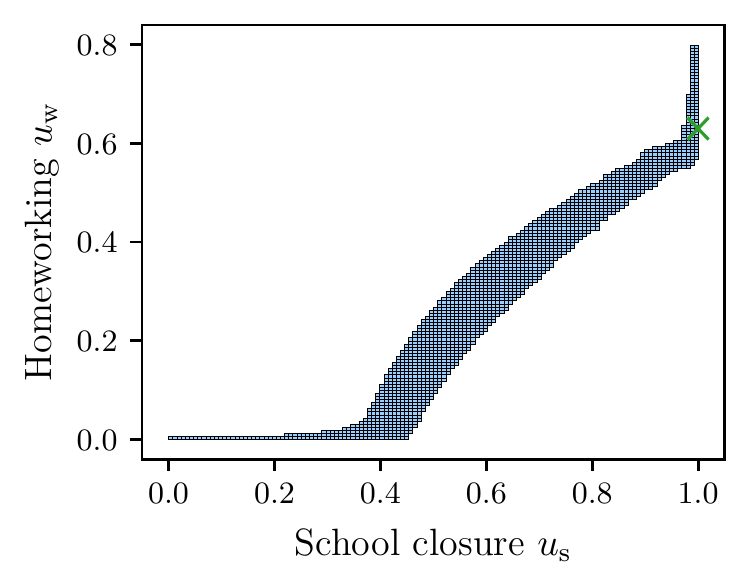}
        \caption{}
    \end{subfigure}
    \begin{subfigure}[t]{0.49\textwidth}
        \includegraphics[width=\textwidth]{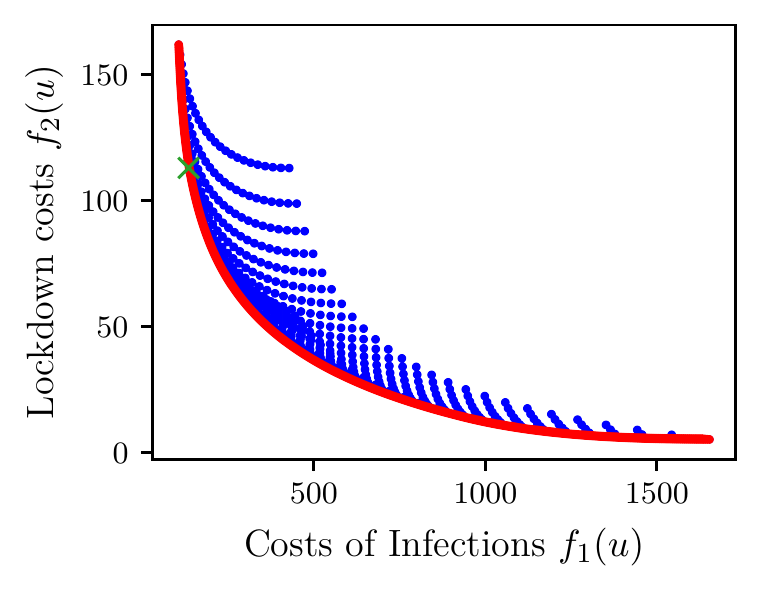}
        \caption{}
    \end{subfigure}
    \caption{(a)~Pareto set covering for feasible decision space $\mathcal{R}$ after 14 subdivision steps obtained with the 6-dimensional reduced model for GERDA. (b)~Pareto front (red) and dominated points (blue). The green cross indicates the optimal control found in~\cite{Niemann2023}.}
    \label{fig:GERDA_pareto_set}
\end{figure}

Since even for smaller geographic regions such as Gangelt, which has roughly 13\ts 000 inhabitants (December 2021), solving the multi-objective optimization problem directly for GERDA is not feasible due to the enormous computational effort and inherent stochasticity of the model. Thus, we construct a surrogate model. Following the arguments in Example~\ref{expl:SIR-ABM}, interpolation between Koopman generators should not be used. In fact, it can be verified easily that gEDMD with control leads to better results (see Table~\ref{tab:GERDA_RMSE}). To train the reduced model, which serves as the surrogate model later on, we sample GERDA (for a fixed world) at 49 points and 225 different control inputs along 7-week trajectories for 24-hourly time steps\new{, i.e., one data point every 24 steps is added to the training data set}. The space of feasible controls $\mathcal{R}$ is discretized using an equidistant grid with $15 \times 15$ grid points. This results in a total of 11\ts 025 training data points, each having 1\ts 000 Monte Carlo simulations to calculate pointwise drift and diffusion estimates using the Kramers--Moyal formulae. We define the augmented state as $x \coloneqq [S_\mathrm{a}, S_\mathrm{c}, I_\mathrm{a}, I_\mathrm{c}, u_\mathrm{s}, u_\mathrm{w}]^\top$, where we distinguish between susceptible and infected adults and children. Note that $R_\mathrm{a}$ and $R_\mathrm{c}$ are omitted due to conservation of $N_\mathrm{a}$ and $N_\mathrm{c}$. The set of basis functions consists of monomials up to degree 4 following the same arguments as in Example~\ref{expl:SIR-ABM}. Figure~\ref{fig:GERDA_reduced} shows the expected aggregated trajectories (dashed) and its standard deviation (shaded) of GERDA as well as the solution of the reduced model (solid) in terms of an ODE for different controls. Table~\ref{tab:GERDA_RMSE} compares the true expected trajectories of GERDA with the reduced model in terms of the root mean squared error for the 6-dimensional state $x$ as well as a further reduced 4-dimensional state $x \coloneqq [S, I, u_\mathrm{s}, u_\mathrm{w}]^\top$. \new{We also include a comparison between the true expected trajectories of GERDA and the generator interpolation method whose use is invalid for the GERDA model. Even on the training data set, state augmentation yields better results.} The 6-dimensional reduced model yields the best approximation. Although all models can be sparsified (e.g., using iterative thresholding~\cite{Brunton2015}), a meaningful relation of the non-zero terms to the classical SIR structure cannot be made. This is, however, also expected as GERDA implicitly contains a network in which agents change locations according to their schedule and interact only with certain agents, i.e., agents present at those locations at the same time. Nevertheless, the reduced model is suitable to solve the multi-objective optimization problem efficiently. Figure~\ref{fig:GERDA_pareto_set} shows the computed Pareto set and front coverings after 14 iterations. In both figures the green cross marks the optimal control $u^\star = [0, 0.63]^\top$ and the corresponding objective values computed in~\cite{Niemann2023}, which was computed without using a surrogate model. We see that this point is included in the Pareto set covering. This finding further highlights the effectiveness of the data-driven surrogate model in approximating the Pareto-optimal set of the ABM. Moreover, solving the multi-objective optimization problem reveals further optimal control inputs belonging to different combinations of weights in~\cite{Niemann2023}.

\begin{table}[]
    \centering
    \caption{Root mean square errors of the trajectories for the reduced models obtained by state augmentation and generator interpolation. The ground truth is a Monte Carlo estimate of GERDA with 1\,000 simulations for each control $u$. The complete parameterization can be found in Table~\ref{tab:parameters}.}
    \label{tab:GERDA_RMSE}
    \begin{tabular}{lcllll}
    \toprule
                            & \multicolumn{1}{l}{}          & \multicolumn{4}{c}{Control $[u_\mathrm{s}, u_\mathrm{w}]^\top$}     \\ \cmidrule{3-6} 
                            & \multicolumn{1}{l}{Dimension} & $[0, 0]^\top$   & $[0, 1]^\top$   & $[1, 0]^\top$  & $[1, 1]^\top$  \\ \midrule
    State augmentation      & 6                             & 0.0325          & 0.0186          & 0.0448         & 0.0023         \\
                            & 4                             & 0.0551          & 0.0325          & 0.0617         & 0.0177         \\ \midrule
    Generator interpolation & 6                             & 0.1869          & 0.1109          & 0.1172         & 0.4122         \\
                            & 4                             & 0.0167          & 0.0347          & 0.0714         & 0.0716         \\ \bottomrule
    \end{tabular}
\end{table}

\begin{table}[]
    \centering
    \caption{Parameters used for the simulations of GERDA.}
    \label{tab:parameters}
    \begin{tabular}{lc}
    \toprule
    Parameter                                              &                      \\ \midrule
    Time $T$ {[}hours{]}                                   & 1176                 \\
    Number of agents $N$                                   & 1045                 \\
    Initially infected $I_\mathrm{a}$                      & 3                    \\
    Initially infected $I_\mathrm{c}$                      & 2                    \\
    World                                                  & Gangelt (reduced)    \\
    General infectivity                                    & 0.175                \\
    General interaction frequency                          & 1                    \\
    Health care system's capacity threshold $I_\text{max}$ & $0.005N$             \\
    Threshold economic impact $u_\mathrm{w}^\text{max}$    & 0.81                 \\ \bottomrule
    \end{tabular}
\end{table}

\begin{rem}
	Note that our primary goal is not to solve multi-objective optimization problems efficiently, but rather to show that surrogate models based on the Koopman generator are well-suited for this purpose. Especially in the case of the voter model, it is more efficient to use the well-known ODE or SDE approximations as surrogate models. However, such approximations are not always known or, in the case of GERDA, might not exist. Thus, constructing surrogate models using the Koopman generator offers a viable approach to solving multi-objective optimization problems involving complex ABMs.
\end{rem}

\section{Conclusions}
\label{sec:conclusion}

The construction of data-driven surrogate models to speed up or even enable the solution of multi-objective optimization problems, where objectives are defined by the dynamical behavior of ABMs, is of great importance as ABMs become more widely used. Especially in the case when limit models are unknown or non-existent, data-driven reduced models and surrogate models allow to solve optimization problems that would otherwise be computationally infeasible due to the very expensive evaluation of the objectives. In this work, we demonstrated how data-driven reduced models based on the Koopman generator enable the efficient solution of multi-objective optimization problems associated with ABMs. For this purpose we constructed surrogate models with varying decision variables, which in this work are controls applied to a non-deterministic dynamical system. We introduced two different methods of constructing these surrogate models. The first one is applicable if the control acts linearly on the system. In this case, linear interpolation between Koopman generators for different control inputs is feasible and can be used to construct surrogate models. This method was illustrated using the voter model, where the goal was to achieve an absolute majority in one opinion with a linear control affecting the second-order transitions. This can be interpreted as repulsion and attraction arguments in an election campaign. The second method overcomes the restriction to linearly acting controls. We showed that a straightforward extension of gEDMD that considers the control as an additional state allows for nonlinear \new{dependence on the} controls. The method was demonstrated first using a simple and analytical model. We then used this method to compute the Pareto set of optimal containment strategies for the GERDA model in a prototypical situation and compared it with the solution computed in~\cite{Niemann2023}. We showed that the solution computed for the single-objective optimization problem is contained in our Pareto set.

In this work, we considered multi-objective optimization problems with $k \ll N$ globally defined objectives, i.e., objectives defined by the collective behavior of the agents, and $n \ll N$ decision variables. An extension is to consider both objectives defined by the individual and collective behavior of agents, which allows for more general problems without considering aggregated trajectories. An important question then is is how to deal with the large computational costs when $k \approx N$ or $n \approx N$ or both. Future research will address these questions.

\subsection*{Acknowledgment}

This research has been funded by the Deutsche Forschungsgemeinschaft (DFG, German Research Foundation) under Germany's Excellence Strategy MATH\texttt{+}: The Berlin Mathematics Research Center (EXC-2046/1, project ID: 390685689).

\subsection*{Declaration of Interests}
The authors declare that they have no known competing financial interests or personal relationships that could have appeared to influence the work reported in this paper.

\bibliographystyle{plain}
% \bibliography{bibfile}

\appendix
\section{\new{Sampling Algorithm}}
\label{sec:appendix_sampling}

For $s=0$, let $\mathcal{B}^0$ denote a collection of finitely many subsets of $\mathcal{R} = [a_1, b_1] \times \dotsc \times [a_n, b_n] \subset \mathbb{R}^n$, $a_i \le y_i \le b_i$ for $i=1,\dots, n$, such that $\bigcup_{B \in \mathcal{B}^0} B = \mathcal{R}$. Then, in each iteration, each box $B \in \mathcal{B}^s$, $s \ge 0$, is first subdivided with respect to one coordinate, resulting in a new collection of boxes $\hat{\mathcal{B}}^{s+1}$. In the second step, a set-wise non-dominance test, which can be carried out heuristically by a finite number of test points, is performed on each box $B \in \hat{\mathcal{B}}^{s+1}$, discarding all boxes containing only dominated points. \new{This process is repeated} with the new collection of boxes $\mathcal{B}^{s+1}$ and the next coordinate. For more details see~\cite{Dellnitz2005, Schuetze2005, Schuetze2013, Peitz2017}. \new{The sampling algorithm is summarized in Algorithm \ref{algo:sampling}.}

\begin{algorithm}
	Let $\mathcal{B}^0$ be a collection of finitely many subsets of $\mathcal{R}$ such that $\bigcup_{B \in \mathcal{B}^0} B = \mathcal{R}$. Then, obtain the new collection $\mathcal{B}^{s+1}$, $s \ge 0$, iteratively from $\mathcal{B}^s$ in two steps:
	\begin{enumerate}[itemsep=0ex, topsep=0.5ex]
		\item Construct a new collection $\hat{\mathcal{B}}^{s+1}$ from $\mathcal{B}^s$ by subdividing each subset $B \in \mathcal{B}^s$\\such that
		\begin{align*}
			\bigcup_{B \in \hat{\mathcal{B}}^{s+1}} B &= \bigcup_{B \in \mathcal{B}^s} B, \\
			\mathrm{diam}(\hat{\mathcal{B}}^{s+1}) &= \theta^{s+1} \mathrm{diam}(\mathcal{B}^s)
		\end{align*}
		for $0 < \theta_\mathrm{min} \le \theta^{s+1} \le \theta_\mathrm{max} < 1$.
		\item Define the new collection $\mathcal{B}^{s+1}$ by
		\begin{equation*}
			\mathcal{B}^{s+1} \coloneqq \left\{ B \in \hat{\mathcal{B}}^{s+1} \mid \nexists \ts \hat{B} \in \hat{\mathcal{B}}^{s+1} \text{ such that } \hat{B} \text{ dominates } B  \right\}.
		\end{equation*}
	\end{enumerate}
	\caption{\cite{Dellnitz2005, Peitz2017}}
	\label{algo:sampling}
\end{algorithm}

\section{\new{Proof of Theorem~\ref{thm:controlaffine}}}
\label{sec:proof}

Given a space of twice continuously differentiable functions and controls $u_1, u_2$, if the dynamics~\eqref{eq:SDE+control} are control-affine, then the Koopman generators are control-affine, i.e.,
\begin{equation*}
	\mathcal{L}_{\alpha_1 u_1 + \alpha_2 u_2} = \mathcal{L}_0 + \alpha_1 \mathcal{A}_{u_1} + \alpha_2 \mathcal{A}_{u_2},
\end{equation*}
where $\mathcal{A}_u = \mathcal{L}_u - \mathcal{L}_0$ and $\alpha_1, \alpha_2 \in \mathbb{R}$.

\begin{proof}
	Consider a control-affine system of the form~\eqref{eq:SDE+control} and let $u$ be any linear combination of controls, i.e., $u = \sum_{i=1}^{d_u} \alpha_i \ts u_i$ for $\alpha_i \in \mathbb{R}$. Given some twice differentiable function $f$, the Koopman generator $\mathcal{L}_u$ depending on $u$ applied to $f$ yields
	\begin{align*}
		\mathcal{L}_u f &= b \cdot \nabla_x f + \frac{1}{2} \ts a \ts \colon \nabla_x^2 f\\
		&= b_0 \cdot \nabla_x f + \sum_{i=1}^{d_u} b_i \ts u_i \cdot \nabla_x f + \frac{1}{2} \ts a_0 \ts \colon \nabla_x^2 f + \frac{1}{2} \sum_{i=1}^{d_u} a_i \ts u_i \ts \colon \nabla_x^2 f,
	\end{align*}
	where $\nabla_x^2$ denotes the Hessian. Defining $\mathcal{A}_u \coloneqq \mathcal{L}_u - \mathcal{L}_0$, we obtain
	\begin{equation*}
		\mathcal{A}_u f = \sum_{i=1}^{d_u} u_i  \left[ b_i \cdot \nabla_x f + \frac{1}{2} \ts a_i \ts \colon \nabla_x^2 f \right].
	\end{equation*}
	The operators $\mathcal{A}_u$ are linear with respect to the control $u$. Moreover, the generators of the Koopman operators are control-affine, that is,
	\begin{equation*}
		\mathcal{L}_u f = \left(\mathcal{L}_0 + \sum_{i=1}^{d_u} \alpha_i \ts \mathcal{A}_{u_i}\right) f. \qedhere
	\end{equation*}
\end{proof}

\end{document}